\documentclass[12pt]{amsart}

\usepackage[utf8]{inputenc}
\usepackage{pxfonts}
\usepackage[margin=1 in]{geometry}

\usepackage[sort&compress,square,numbers]{natbib}
\usepackage{graphicx,amssymb,amstext,amsmath,wrapfig,url,bm}
\usepackage[colorlinks]{hyperref}
\usepackage[hang]{subfigure}

\renewcommand{\vec}[1]{\mathbf{#1}}

\usepackage{amsmath,amsthm}

\usepackage{theomac}

\newcommand{\lemlab}[1]{\label{lemma:#1}}
\newcommand{\theolab}[1]{\label{theo:#1}}

\newcommand{\eqlab}[1]{\label{eq:#1}}
\newcommand{\corlab}[1]{\label{cor:#1}}

\newcommand{\seclab}[1]{\label{section:#1}}

\newcommand{\proplab}[1]{\label{prop:#1}}

\newcommand{\lemref}[1]{Lemma \ref{lemma:#1}}

\newcommand{\theoref}[1]{Theorem \ref{theo:#1}}

\newcommand{\corref}[1]{Corollary \ref{cor:#1}}

\renewcommand{\eqref}[1]{(\ref{eq:#1})}
\newcommand{\secref}[1]{Section \ref{section:#1}}

\newcommand{\propref}[1]{Proposition \ref{prop:#1}}

\newtheoremWithMacro{theorem}{Theorem}

\newtheorem{lemma}{Lemma}
\newtheoremWithMacro{prop}[lemma]{Proposition}
\newtheorem{cor}[lemma]{Corollary}

\begin{document}

\title{Natural Realizations of Sparsity Matroids}

\author{Ileana Streinu}
\address{Computer Science Department, Smith College \\ Northampton, MA}
\email{streinu@cs.smith.edu,istreinu@smith.edu}
\urladdr{http://cs.smith.edu/~streinu}
\author{Louis Theran}
\address{Mathematics Department \\ 
Temple University, Philadelphia, PA}
\email{theran@temple.edu}
\urladdr{http://math.temple.edu/~theran/}

\begin{abstract}
A hypergraph $G$ with $n$ vertices and $m$ hyperedges with $d$ endpoints each 
is $(k,\ell)$-sparse if for all sub-hypergraphs $G'$ on $n'$ vertices and $m'$ edges, $m'\le kn'-\ell$.  For 
integers $k$ and $\ell$ satisfying $0\le \ell\le dk-1$, this is known to be a linearly representable matroidal
family.

Motivated by problems in rigidity theory, we give a new linear representation theorem for the $(k,\ell)$-sparse hypergraphs that is \emph{natural}; i.e., the representing matrix captures the vertex-edge incidence structure 
of the underlying hypergraph $G$.
\end{abstract}

\maketitle

%%%%%%%%%%%%%%%%%%%%%%%%%%%%%%%%%%%%%%%%%%%%%%%%%%%%%%%%%%%%%%%%%%%%%
%%%%%%%%%%%%%%%%%%%%%%%%%%%%%%%%%%%%%%%%%%%%%%%%%%%%%%%%%%%%%%%%%%%%%
\section{Introduction}
Let $G$ be a \emph{$d$-uniform hypergraph}; i.e., $G=(V,E)$, where $V$ is a 
finite set of $n$ \emph{vertices} and $E$ is a multi-set of $m$ \emph{hyperedges}, 
which each have $d$ distinct \emph{endpoints}.  We define $G$ to be 
\emph{$(k,\ell)$-sparse} if, for fixed integer parameters $k$ and $\ell$, 
any sub-hypergraph $G'$ of $G$ on $n'$ vertices 
and $m'$ hyperedges satisfies the relation $m'\le kn'-\ell$; 
if, in addition $m=kn-\ell$, then $G$ is \emph{$(k,\ell)$-tight}.

For a fixed $n$, and integer parameters $k,\ell$, and $d$ satisfying $0\le \ell\le dk-1$,
the family of $(k,\ell)$-tight $d$-uniform hypergraphs on $n$ vertices 
form the bases of a matroid \cite{whiteley:matroids}, which we define to be the 
\emph{$(k,\ell)$-sparsity-matroid}.  The topic of this paper is linear 
representations of the $(k,\ell)$-sparsity-matroids with a specific 
form.  

\subsection*{Main Theorem.} Our main result is the following.  Detailed definitions of 
$(k,\ell)$-sparse hypergraphs are given in \secref{sparse}; detailed
definitions of linear representations are given in \secref{natural}.

\begin{theorem}[\main][\textbf{Natural Realizations}]\theolab{main}
Let $k$, $\ell$, and $d$ be integer parameters satisfying the  
inequality $0\le \ell\le kd-1$.  Then, for sufficiently large $n$,
the $(k,\ell)$-sparsity-matroid of $d$-uniform 
hypergraphs on $n$ vertices is representable 
by a matrix $\vec M$ with:
\begin{itemize}
	\item Real entries
	\item $k$ columns corresponding to each vertex (for a total of $kn$)
	\item One row for each hyperedge $e$
	\item In the row corresponding to each edge $e$, the only non-zero entries 
	appear in columns corresponding to endpoints of $e$
\end{itemize}
\end{theorem}

\subsection*{Novelty.}  
As a comparison, standard matroidal constructions imply that there is a linear 
representation that is $m\times kn$ for all the allowed values of $k$, $\ell$ and $d$.
For $d=2$, $\ell\le k$, the $(k,\ell)$-sparsity-matroid is characterized as the matroid union 
of $\ell$ copies of the standard graphic matroid and $(k-\ell)$ copies of the bicycle matroid,
so the desired representation follows from the Matroid Union Theorem \cite[Section 7.6]{constructions} 
for linearly representable matroids.

\theoref{main}, in contrast, applies to \emph{the entire matroidal range of parameters $k$, $\ell$, and $d$}.
In particular, it applies in the so-called \emph{upper range} in which $\ell> k$.  In the upper range, 
no reduction to matroid unions are known, so proofs based on the Matroid Union Theorem do not apply.

\subsection*{Motivation.}  Our motivation for this work comes from \emph{rigidity theory},
which is the study of structures defined by geometric constraints.  Examples include:  \emph{bar-joint frameworks},
which are structures made of \emph{fixed-length bars} connected by \emph{universal joints}, with full rotational freedom; and \emph{body-bar frameworks}, which are made of \emph{rigid bodies} connected by fixed length bars attached to universal joints.  A framework is \emph{rigid} if the only allowed continuous motions 
that preserve the lengths and connectivity of the bars are rigid motions of Euclidean space.

In both cases, the formal description of the framework is given in two parts: 
a \emph{graph} $G$, defining the combinatorics of the framework; \emph{geometric data}, 
specifying the lengths of the bars, and their attachment points on the bodies.  Rigidity 
is a difficult property to establish in all cases, the with best known algorithms 
relying on exponential-time Gröbner basis computations.  However, for \emph{generic}
geometric data (and almost all lengths are generic, see \cite{sliders} for a detailed discussion), 
rigidity properties can be determined from the combinatorics of the framework alone, as shown by the 
following two landmark theorems:

\begin{theorem}[\laman][\textbf{Maxwell-Laman Theorem: Generic planar bar-joint rigidity \cite{laman,maxwell}}]\theolab{laman}
A generic bar-joint framework in $\mathbb{R}^2$ is minimally rigid if and only if 
its underlying graph $G$ is $(2,3)$-tight.
\end{theorem}

\begin{theorem}[\tay][\textbf{Tay's Theorem: Generic body-bar rigidity \cite{tay}}]\theolab{tay}
A generic body-bar framework in $\mathbb{R}^d$ is minimally rigid if and only if its 
underlying graph $G$ is $(\binom{d+1}{2},\binom{d+1}{2})$-tight.
\end{theorem}

All known proofs of theorems such as \ref{theo:laman} and \ref{theo:tay} proceed via a linearization 
of the problem called \emph{infinitesimal rigidity}.  The key step in all of these proofs is to prove that 
a specific matrix, called the \emph{rigidity matrix}, which arises as the differential of the equations 
for the length constraints, is, generically, a linear representation of some $(k,\ell)$-sparsity matroid.

The rigidity matrices arising in Theorems \ref{theo:laman} and \ref{theo:tay} are specializations 
of our natural realizations: they have the same pattern of zero and non-zero entries.  The present 
work arises out of a project to understand \emph{``rigidity from the combinatorics up''} by 
studying \emph{$(k,\ell)$-sparse graphs} and their generalizations.  Our main \theoref{main}
and the implied natural realizations occupy an intermediate position in between the rigidity 
theorems and the combinatorial matroids of $(k,\ell)$-sparse graphs.  The natural realizations 
presented here may be useful as building blocks for a new, more general class 
of rigidity theorems in the line of \ref{theo:laman} and \ref{theo:tay}.

\subsection*{Related work: $(k,\ell)$-sparse graphs.}  Graphs and hypergraphs 
defined by hereditary sparsity counts first appeared as an example of matroidal 
families in the the work of Lorea \cite{lorea}.  Whiteley, as part of a project with Neil White, 
reported in  \cite[Appendix]{whiteley:matroids}, studied them from the rigidity perspective.  
Michael Albertson and Ruth Haas \cite{albertson:haas} studied $(k,\ell)$-sparse graphs from an extremal 
perspective as an instance of graphs characterized by ``bounding functions.''

This paper derives more directly from the sequence of papers by Ileana Streinu and 
her collaborators: \cite{pebblegame} develops the structural and algorithmic theory
of $(k,\ell)$-sparse graphs; \cite{hypergraphs} extends the results of \cite{pebblegame}
to hypergraphs; \cite{maps,colors} give characterizations in terms of decompositions 
into trees and ``map-graphs''; \cite{graded} extends the sparsity concept to allow
different counts for different types of edges.

\subsection*{Related work: matroid representations.}  For the specific 
parameter values $d=2$, $\ell\le k$, natural realizations of the type 
presented in \theoref{main} may be deduced from the Matroid Union Theorem \cite[Section 7.6]{constructions};
this was done by Whiteley \cite{whiteley:union-matroids}, where the realizations for $d=2$, $\ell=l$ go 
by the name ``$k$-frame.''  In addition, White and Whiteley \cite[Appendix]{whiteley:matroids} 
have shown, using a geometric construction involving picking projective flats in general 
position and then the Higgs Lift \cite[Section 7.5]{constructions} 
that all $(k,\ell)$-sparsity matroids for graphs and hypergraphs are linearly representable.
%XXX

Whiteley \cite{Whiteley:1989p992}
proved a very similar result for the special case of $k=1$; he also 
gave representations for a related class of matroids on bipartite incidence 
graphs\footnote{These matroids 
have also appeared in the Ph.D. thesis of Audrey Lee-St. John \cite{audreythesis} 
under the name ``mixed sparsity.''}.

All known rigidity representation theorems \cite{laman,kt,sliders,tay,whiteley:union-matroids} 
provide natural realizations for the specific 
sparsity parameters involved.  However, all these give \emph{more specialized} representations,
arising from geometric considerations, with more specialized proofs.  All the arguments having a 
matroidal flavor seem to rely, in one way or another, on the Matroid Union Theorem, or the 
explicit determinantal formulas used to prove it.

\subsection*{Related work: rigidity theory.}  Lovász and Yemini \cite{ly} introduced the matroidal 
perspective to rigidity theory with their proof of the Maxwell-Laman \theoref{laman} based on an
explicit computation of the rank function of the $(2,3)$-sparsity matroid that uses its special 
relationship with the union of graphic matroids.  Whiteley \cite{whiteley:union-matroids} gives a very elegant proof 
of Tay's \theoref{tay} \cite{tay} using the Matroid Union Theorem and geometric 
observations specific to the body-bar setting.  White and Whiteley \cite{white:whiteley} analyzed the 
minors of $k$-frames of \cite{whiteley:union-matroids} in detail, 
describing ``pure conditions'' that determine the rigidity behavior of body-bar frameworks.

In both \cite{ly,whiteley:union-matroids}, as well as in more recently proven 
Maxwell-Laman-type theorems of Katoh and Tanigawa \cite{kt} 
and the authors' \cite{sliders}, the connection between $(k,\ell)$-sparsity and 
\emph{sparsity-certifying decompositions} \cite{colors} of the minimally rigid family of 
graphs appears in an essential way.  In contrast, here we only need to employ 
sparsity itself, yielding a much more general family of realizations.  The price for 
this added generality is that we cannot immediately deduce rigidity results directly 
from \theoref{main}.

\subsection*{Organization.}   \secref{sparse} introduces $(k,\ell)$-sparse
hypergraphs and gives the necessary structural properties.  \secref{natural} gives the required background 
in linear representability of matroids and then the proof of \theoref{main}.  In \secref{extensions},
we describe two extensions of \theoref{main}: to non-uniform $(k,\ell)$-sparse hypergraphs and 
to $(k,\bm{\ell})$-graded-sparse hypergraphs.  We conclude in \secref{conclusions} with some remarks 
on the relationship between natural realizations and rigidity.

\subsection*{Notations.}  A \emph{hypergraph} $G=(V,E)$ is defined by a finite set $V$ of \emph{vertices} and a 
multi-set $E$ of \emph{hyperedges}, which are subsets of $V$; if $e \in E(G)$ is an edge and $v\in e$ is a 
vertex, then we call $v$ an \emph{endpoint} of the edge $e$.  A hypergraph $G$ is defined to be $d$-uniform if
all the edges have $d$ endpoints.  Sub-hypergraphs are typically denoted as $G'$ with $n'$ vertices and $m'$ edges; whether they are 
vertex- or hyperedge-induced will be explicitly states.  For $d$-uniform hypergraphs, we use the notation $e_1,e_2,\ldots,e_d$ for the 
$d$ endpoints of a hyperedge $e\in E(G)$.

Matrices $\vec M$ are denoted by bold capital letters, vectors $\vec v$ by bold lowercase letters.  The 
rows of a matrix $\vec M$ are denoted by $\vec m_i$.

The letters $k$, $\ell$, and $d$ denote sparsity parameters.

\subsection*{Dedication.} This paper is dedicated to the memory of Michael Albertson.

\section{The $(k,\ell)$-sparsity matroid}\seclab{sparse}
Let $(k,\ell,d)$ be a triple of non-negative integers such that $0\le \ell\le dk-1$; we define such 
a triple as giving \emph{matroidal sparsity parameters} 
(this definition is justified below in \propref{matroid}).  
A $d$-uniform hypergraph $G=(V,E)$ with $n$ vertices and $m$ hyperedges is\emph{ $(k,\ell)$-sparse} if, 
for all subsets $V'\subset V$ of $n'$ vertices, the subgraph induced by $V'$ has $m'$ edges with 
$m'\le kn'-\ell$.  If, in addition, $m=kn-\ell$, $G$ is \emph{$(k,\ell)$-tight}.  
For brevity, we call $(k,\ell)$-tight $d$-uniform hypergraphs $(k,\ell,d)$-graphs.

The starting points for the results of this paper is the matroidal property of 
$(k,\ell,d)$-graphs.  We define $K_{n,d}^{dk-\ell}$ to be the complete 
$d$-uniform hypergraph on $n$ vertices with $dk-\ell$ copies of each hyperedge.
\begin{prop}[\matroidal][\cite{lorea,whiteley:matroids,hypergraphs}]\proplab{matroid}
Let $d$, $k$ and $\ell$ be non-negative integers satisfying $\ell\in [0,dk-1]$.  Then the family 
of $(k,\ell,d)$-graphs on $n$ vertices forms the bases of a matroid on the edges of $K_{n,d}^{dk-\ell}$, 
for a sufficiently large $n$, depending on $k$, $\ell$, and $d$.
\end{prop}
We define the matroid appearing in \propref{matroid} to be the \emph{$(k,\ell,d)$-sparsity-matroid}.

From now on, the parameters $k$, $\ell$ and $d$ are always matroidal sparsity parameters 
and $n$ is assumed to be large enough that \propref{matroid} holds.

The other fact we need is the following lemma from \cite{hypergraphs} characterizing the special case
of $(k,0,d)$-graphs.  We define an \emph{orientation} of a hypergraph to be an assignment of a \emph{tail}
to each hyperedge by selecting one of its endpoints (unlike in the graph setting, there is no uniquely defined 
head).
\begin{lemma}[\cite{hypergraphs}]\lemlab{orientation}
Let $G$ be a $d$-uniform hypergraph with $n$ vertices and $m=kn$ hyperedges.  Then $G$ is a
$(k,0,d)$-graph if and only if there is an orientation such that each vertex is the tail of 
exactly $k$ hyperedges.
\end{lemma}

\section{Natural Realizations}\seclab{natural}
In this section we prove our main theorem:
\main

\subsection*{Roadmap.}  This section is structured as follows.  We begin by defining generic matrices and 
then introduce the required background in linear representation of matroids.  The proof of \theoref{main}
then proceeds by starting with the special case of $(k,0)$-sparse hypergraphs and then reducing to 
it via a general construction.

\subsection*{The generic rank of a matrix.}  A \emph{generic matrix}
has as its non-zero entries \emph{generic variables}, or formal polynomials over 
$\mathbb{R}$ or $\mathbb{C}$ in generic variables.  
Its \emph{generic rank} 
is given by the largest number $r$
for which $\vec M$ has an $r\times r$ matrix minor with a determinant that is formally non-zero.

Let $\vec M$ be a generic matrix in $m$ generic variables $x_1,\ldots, x_m$, and 
let $\vec v=(v_i)\in \mathbb{R}^m$ (or $\mathbb{C}^m$).  We define a \emph{realization of $\vec M$} to be the 
matrix obtained by replacing the variable $x_i$ with the corresponding number $v_i$.  A 
vector $\vec v$ is defined to be a \emph{generic point} if the rank of the associated realization 
is equal to the  generic rank of $\vec M$; otherwise $\vec v$ is defined to be a \emph{non-generic} point.

We will make extensive use of the following well-known facts from algebraic geomety (see, e.g., \cite{iva}):
\begin{itemize}
	\item The rank of a generic matrix $\vec M$ in $m$ variables 
	is equal to the maximum over $\vec v\in \mathbb{R}^m$ ($\mathbb{C}^m$) of the rank of all 
	realizations.
	\item The set of non-generic points of a generic matrix $\vec M$ is an 
	algebraic subset of $\mathbb{R}^m$ ($\mathbb{C}^m$).
	\item The rank of a generic matrix $\vec M$ in $m$ variables is at least 
	as large as the rank of any specific realization; i.e., generic rank 
	can be established by a single example.
\end{itemize}

\subsection*{Generic representations of matroids.}
Let $\mathcal{M}$ be a matroid on ground set $E$.  We define a generic matrix $\vec M$ to be 
a \emph{generic representation of $\mathcal{M}$} if:
\begin{itemize}
	\item There is a bijection between the rows of $\vec M$ and the ground set $E$.
	\item A subset of rows of $\vec M$  attains the rank of the matrix $\vec M$ 
	if and only if the corresponding subset of $E$ is a basis of $\mathcal{M}$.
\end{itemize}

\subsection*{Natural realizations for $(k,0,d)$-graphs.}  Fix matroidal parameters $k$, $\ell=0$ and $d$, and 
let $G$ be a $d$-uniform hypergraph on $n$ vertices and $m$ hyperedges.  For a hyperedge $e\in E(G)$ with endpoints 
$e_i$, $i\in [1,d]$, define the vector $\vec a_{e_i} = (a^j_{e_i})_{j\in [1,k]}$
to have as its entries $k$ generic variables for each of the $d$ endpoints of $e$.  

Next, we define the generic matrix $\vec M_{k,0,d}(G)$ to have $m$ rows, indexed by the hyperedges of $G$,
and $kn$ columns, indexed by the vertices of $G$, with $k$ columns for each vertex.  The filling pattern
of $\vec M_{k,0,d}$ is given as follows:
\begin{itemize}
	\item If a vertex $i\in V(G)$ is an endpoint of an edge $e$, then the $k$ entries associated with $i$
	in the row indexed by $e$ are given by the vector $\vec a_{e_i}$.
	\item All other entries are zero.
\end{itemize}
For example, if $G$ is a $3$-uniform hypergraph, the matrix $\vec M_{k,0,3}(G)$ has the following pattern:
\[  \bordermatrix{                   &              &     e_1    &             &     e_2 &            &           e_3                               &  \cr
					   & \cdots & \cdots               & \cdots & \cdots               & \cdots & \cdots                                      & \cdots \cr
					e   &  0\cdots0 & a^1_{e_1} \cdots a^k_{e_1} & 0\cdots0 & a^1_{e_2}\cdots a^k_{e_2} & 0\cdots 0& a^1_{e_3} \cdots a^k_{e_3} & 0 \cdots 0\cr
					   & \cdots & \cdots                & \cdots & \cdots               & \cdots & \cdots                                     & \cdots  }.
\]

The following lemma is a consequence of the Matroid Union Theorem and a 
representation result for the $(1,0,d)$-sparsity-matroid due to 
Edmonds \cite{edmonds}.\footnote{Whiteley \cite[Prop. 2.4]{Whiteley:1989p992} 
reproduces Edmonds's proof; here, even in the $(1,0,d)$-case we 
go along different lines.}  We give a more direct proof for completeness.
\begin{lemma}\lemlab{k0dtight}
Let $G$ be a $d$-uniform hypergraph on $n$ vertices and $m=kn$ edges.  Then 
$\vec M_{k,0,d}(G)$ has generic rank $kn$ if and only if $G$ is a $(k,0,d)$-graph.
\end{lemma}
\begin{proof}
First, we suppose that $G$ is a $(k,0,d)$-graph.  By \lemref{orientation}, there is 
an assignment of a distinct tail to each edge such that each vertex is the tail of 
exactly $k$ edges.  Fix such an orientation, giving a natural association of $k$ edges
to each vertex.  Now specialize the matrix $\vec M_{k,0,d}(G)$ as follows:
\begin{itemize}
	\item Let $i\in V(G)$ be a vertex that is the tail of edges $e_{i_1},e_{i_2},\ldots,e_{i_k}$.
	\item In row $e_{i_j}$, set the variable $a^j_{e_{i_j}}$ to 1 and all other entries to zero.
\end{itemize}
Because each edge has exactly one tail, this process defines a setting for the entries of $\vec M_{k,0,d}(G)$
with no ambiguity.  Moreover, after rearranging the rows and columns, this setting of the 
entries turns $\vec M_{k,0,d}(G)$ into the identity matrix, so this example shows its rank generic is $kn$.

In the other direction, we suppose that $G$ is not a $(k,0,d)$-graph.  Since $G$ has $kn$ edges,
it is not $(k,0)$-sparse, so some subgraph $G'$ spanning $n'$ vertices induces at least $kn'+1$
edges.  Arranging the edges and vertices of $G'$ into the top-left corner of $\vec M_{k,0,d}(G)$,
we see that $G'$ induces a submatrix with at least $kn'+1$ rows and only $kn'$ columns that are not 
entirely zero.  It follows that $\vec M_{k,0,d}(G)$ must be, generically, rank deficient.
\end{proof}

\begin{cor}\corlab{k0rep}
The matrix $\vec M_{k,0,d}(K_{n,d}^{dk})$ is a generic representation for the $(k,0,d)$-sparsity matroid.	
\end{cor}
\begin{proof}
\lemref{k0dtight} shows that a $kn\times kn$ matrix minor is generically non-zero if and only if the 
set of rows it induces corresponds to a $(k,0,d)$-graph, so the bases of $\vec M_{k,0,d}(K_{n,d}^{dk})$
are in bijective correspondence with $(k,0,d)$-graphs. 
\end{proof}

\begin{cor}\corlab{k0gen}
Let $G$ be a $(k,\ell)$-sparse $d$-uniform hypergraph with $m$ hyperedges.
The set of $\vec v\in \mathbb{R}^{dkm}$ such that the associated realization of 
$\vec M_{k,0,d}(G)$ has full rank is the open, dense complement of an algebraic subset of $\mathbb{R}^{dkm}$.
\end{cor}
\begin{proof}
\corref{k0rep} implies that the rank drops only when $\vec v$ is a common zero of all the $m\times m$
minors of $\vec M_{k,0,d}(G)$, which is a polynomial condition.
\end{proof}

\subsection*{The natural representation matrix $\vec M_{k,\ell,d}(G)$.}  
Fix matroidal sparsity parameters $k$, $\ell$, and $d$, and let
$G$ be a $d$-uniform hypergraph. Let $\vec U$ be an $kn\times\ell$ matrix with generic entries.  
We define the matrix $\vec M_{k,\ell,d}(G)$ to be a generic matrix that is a formal solution to the equation \eqref{cut} below, with the
entries of $\vec U$ fixed and the entries of $\vec M_{k,0,d}(G)$ as the variables:
\begin{equation}\eqlab{cut}
	\vec M_{k,0,d}(G)\vec U = 0
\end{equation}
We note that the process of solving \eqref{cut} does not change the location of zero and non-zero entries in $	\vec M_{k,0,d}(G)$, preserving the 
naturalness property required by \theoref{main}.

With this definition, we can restate \theoref{main} as follows: 
the matrix $\vec M_{k,\ell,d}(K_{n,d}^{dk-\ell})$ is a generic representation of the $(k,\ell,d)$-sparsity 
matroid.

\subsection*{Main lemmas}
The next two lemmas give the heart of the proof of \theoref{main}.  The first says that 
if $G$ is not $(k,\ell)$-sparse, then $\vec M_{k,\ell,d}(G)$ has a row dependency.

\begin{lemma}\lemlab{rankdrop}
Let $k$, $\ell$, and $d$ be matroidal parameters and be $G$ a  $d$-uniform hypergraph with $m=kn-\ell$.  
If $G$ is not $(k,\ell)$-sparse, then $\vec M_{k,\ell,d}(G)$ is not generically full rank.
\end{lemma}
\begin{proof}
Since $G$ is not $(k,\ell)$-sparse, it must have some vertex-induced 
subgraph $G'$ on $n'$ vertices and $m'>kn'-\ell$ edges.  
The  sub-matrix of $\vec M_{k,\ell,d}(G)$ induced by the edges of $G'$ has at least $kn'-\ell+1$ 
rows and only $kn'$ columns that are not all zero, so it must have a row dependency, since, by 
definition, the kernel of such a sub-matrix has dimension at least $\ell$.
\end{proof}

The following is the key lemma.  It says that, generically, the dependencies of the type
described by \lemref{rankdrop} are the only ones.
\begin{lemma}\lemlab{ranknodrop}
Let $k$, $\ell$, and $d$ be matroidal parameters and be $G$ a  $d$-uniform hypergraph with $m=kn-\ell$.  
If $G$ is $(k,\ell)$-sparse, i.e., it is a $(k,\ell,d)$-graph, 
then $\vec M_{k,\ell,d}(G)$ is  generically full rank.
\end{lemma}
\begin{proof}
We prove this by constructing an example, from which the generic statement follows.
From \corref{k0rep} and \corref{k0gen}, we may select values for the variables 
$a^j_{e_i}$ in the generic matrix $\vec M_{k,0,d}(G)$ so that the resulting 
realization $\vec M$ of $\vec M_{k,0,d}(G)$ is full rank.  

Denote by $\vec m_e$, for $e\in E(G)$, the rows of $\vec M$.  
Define the subspace $W_G$ of $\mathbb{R}^{kn}$ to be the linear span of the $\vec m_e$.  
For each vertex-induced subgraph $G'$ on $n'$ vertices of $G$ define $W_{G'}$ to be the linear 
span of $\{\vec m_e : e\in E(G')\}$; $W_{G'}$ is a subspace of $\mathbb{R}^{kn}$, and, 
because the $\vec m_e$ span exactly $kn'$ non-zero columns in $\vec M$, 
it has a natural identification as a subspace of $\mathbb{R}^{kn'}$.

We will show that there is a subspace $U$ of $\mathbb{R}^{kn}$ such that $W_{G}\cap U^\perp$
has dimension $kn-\ell$; taking the matrix $\vec U$ to be a basis of $U$ and then solving 
$\vec m_e\vec U = 0$ for each row of $\vec M$ gives a solution to \eqref{cut} with full 
rank.  This proves the lemma, since the resulting matrix will have as its rows a basis
for $W_{G}\cap U^\perp$, which has dimension $kn-\ell$.

Now let $U$ be an $\ell$-dimensional subspace of $\mathbb{R}^{kn}$ with basis given by the columns of
the $kn\times \ell$ matrix $\vec U$.  For each vertex-induced subgraph $G'$ of $G$ on  $n'$ vertices, 
associate the corresponding $kn'$ rows of $\vec U$ to determine a subspace $U_{G'}$.

Let $G'$ be a vertex-induced subgraph of $G$ on $n'$ vertices 
and consider the subspace $W_{G'}$.  
Since $\operatorname{dim} W_{G'} = \operatorname{dim} (W_{G'}\cap U_{G'}) + 
\operatorname{dim} (W_{G'}\cap U^\perp_{G'})$, 
if $\operatorname{dim} (W_{G'}\cap U^\perp_{G'}) < \operatorname{dim} W_{G'}$, then$ W_{G'}\cap U_{G'}$ is at least one-dimensional.

Here is the key to the proof (and where the combinatorial assumption of $(k,\ell)$-sparsity
enters in a fundamental way): by the $(k,\ell)$-sparsity of $G$, 
the dimension of $W_{G'}$ is at most $kn'-\ell$.  Since $U_{G'}$ is only (at most) an $\ell$-dimensional 
subspace of $\mathbb{R}^{kn'}$, this only happens if the bases of $W_{G'}$ and $U_{G'}$ satisfy
a polynomial relation.  Since there are only finitely many subgraphs, this gives a finite polynomial 
condition specifying which $U$ are disallowed, completing the proof.
\end{proof}
	
\subsection*{Proof of the Main \theoref{main}} With the two key Lemmas \ref{lemma:rankdrop} and \ref{lemma:ranknodrop}, 
the proof of \theoref{main} is very similar to that of \corref{k0rep}.  
We form the generic matrix $\vec M_{k,\ell,d}(K_{n,d}^{dk-\ell})$.  \lemref{rankdrop} and 
\lemref{ranknodrop} imply that a set of rows forms a basis if and only if the corresponding hypergraph $G$
is a $(k,\ell,d)$-graph.
\hfill $\qed$

\section{Extensions: non-uniform hypergraphs and graded sparsity}\seclab{extensions}
In this section, we extend \theoref{main} in two directions: to $(k,\ell)$-sparse hypergraphs
that are not $d$-uniform; to $(k,\bm{\ell})$-graded sparse hypergraphs.

\subsection*{Non-uniform hypergraphs.}  The theory of $(k,\ell)$-sparsity we developed in 
\cite{hypergraphs}, does not require that a hypergraph $G$ be $d$-uniform.  All the definitions
are similar, except we require only that if $\ell\ge (d-1)k$, then each hyperedge have at least $d$
endpoints.  The ground set of the corresponding sparsity matroid now is the more complicated
hypergraph on $n$ vertices with $ik-\ell$ copies of each hyperedge with $i$ endpoints for $i\ge d$.

The combinatorial properties enumerated in \secref{sparse} all hold in the non-uniform setting, 
and the proofs in \secref{natural} all go through verbatim, with slightly more complicated notation, yielding:
\begin{theorem}[\nonuniform][\textbf{Natural Realizations: non-uniform version}]\theolab{nonuniform}
Let $k$, $\ell$, be integer parameters satisfying the  
inequality $0\le \ell\le kd-1$.  Then, for sufficiently large $n$,
the $(k,\ell)$-sparsity-matroid of non-uniform hypergraphs on $n$ vertices is representable 
by a matrix $\vec M$ with:
\begin{itemize}
	\item Real entries
	\item $k$ columns corresponding to each vertex (for a total of $kn$)
	\item One row for each hyperedge $e$
	\item In the row corresponding to each edge $e$, the only non-zero entries 
	appear in columns corresponding to endpoints of $e$
\end{itemize}
\end{theorem}

\subsection*{Graded-sparsity.}  In \cite{graded}, we developed an extension of $(k,\ell)$-sparsity
called $(k,\bm{\ell})$-graded-sparsity.  Graded-sparsity is the generalization of the 
sparsity counts appearing in our work on slider-pinning rigidity \cite{sliders}.

Define the hypergraph $K^+_{n,k}$, to be complete hypergraph on $n$ vertices, where hyperedges 
with  $d$ endpoints have multiplicity $dk$. A {\em grading} $(E_1,E_2,\ldots,E_s)$ of $K^+_n$ 
is a strictly decreasing sequence of sets of edges $E(K^+_n)=E_1\supsetneq E_2\supsetneq \cdots\supsetneq E_s$.
Now fix a grading on $K^+_n$ and let $G=(V,E)$ be a hypergraph.
Define $G_{\ge i}$ as the subgraph of $G$ induced by $E\cap E_j$. % $E\cap \left(\cup_{j\ge i}E_i\right)$. 
Let $\bm{\ell}$
be a vector of $s$ non-negative integers.  We say that $G$ is $(k,\bm{\ell})$-graded sparse if
$G_{\ge i}$ is $(k,\bm{\ell}_i)$-sparse for every $i$; $G$ is {\bf $(k,\bm{\ell})$-graded tight}
if, in addition, it is $(k,\bm{\ell}_1)$-tight.

The main combinatorial result of \cite{graded} is that $(k,\bm{\ell})$-graded-sparse hypergraphs form the 
bases of a matroid, which we define to be the $(k,\bm{\ell})$-graded-sparsity matroid.

\begin{theorem}[\graded][\textbf{Natural Realizations: graded-sparsity}]\theolab{graded}
Fix a grading of $K_{n,k}^+$ and let $k$ and $\bm{\ell}$ be graded-sparsity parameters.
Then, for sufficiently large $n$,
the $(k,\bm{\ell})$-sparsity-matroid s on $n$ vertices is representable 
by a matrix $\vec M$ with:
\begin{itemize}
	\item Real entries
	\item $k$ columns corresponding to each vertex (for a total of $kn$)
	\item One row for each hyperedge $e$
	\item In the row corresponding to each edge $e$, the only non-zero entries 
	appear in columns corresponding to endpoints of $e$
\end{itemize}
\end{theorem}

Because of the presence of the grading, we need to modify the proof of \theoref{main} to account for 
it.  The formal matrix $\vec M_{k,0,+}(K^+_{n,k})$ is defined analogously to 
$\vec M_{k,0,d}(K_{n,d}^{dk-\ell})$, except we sort the rows by the grading.  The counterpart to 
\eqref{cut} then becomes the system:
\begin{equation}\eqlab{gradedcut}
\vec M_{k,0,d}(E_{\ge i}) \vec U_i = 0, \,\, i=1,2,\ldots, s
\end{equation}
where $V_1$ is $kn\times \ell_1$, and each successive $\vec U_i$ is $\vec U_i$ with $\ell_i$
additional columns.

With this setup, the proof of \theoref{main} goes through with appropriate notational changes.

\section{Conclusions and remarks on rigidity}\seclab{conclusions}
We provided linear representations for the matroidal families of $(k,\ell)$-sparse hypergraphs and $(k,\bm{\ell})$-graded-sparse hypergraphs that are \emph{natural}, in the sense that the representing
matrices capture the vertex-edge incidence pattern.  This family of representations, which 
extends to the \emph{entire matroidal range} of sparsity parameters, may be useful as a building 
block for ``Maxwell-Laman-type'' rigidity theorems.  We conclude with a brief discussion of why
one \emph{cannot} conclude rigidity theorems such as \theoref{laman} and \theoref{tay}
directly from \theoref{main}.

The proof of the critical \lemref{ranknodrop} is very general, since it has to work for the entire
range of sparsity parameters.  What it guarantees is that the entries of $\vec M_{k,\ell,d}(G)$
are some polynomials, but not what these polynomials are.  For rigidity applications, 
\emph{specific} polynomials are forced by the geometry, which would require more control over
the matrix $\vec U$ appearing in Equation \eqref{cut} than the proof technique here allows.

For example, in the planar bar-joint rigidity case the ``trivial infinitesimal motions'' can be given the basis:
\begin{itemize}
	\item $(1,0,1,0\ldots,1,0)$ and $(0,1,0,1,\ldots,0,1)$, representing infinitesimal translation
	\item $(-y_1,-y_2,\ldots,-y_n,x_1,x_2,\ldots,x_n)$, representing infinitesimal rotation around the origin
\end{itemize}

It is important to note that \theoref{main} \emph{cannot} simply be applied with this collection as the 
columns of $\vec U$ to conclude the Maxwell-Laman Theorem \ref{theo:laman}.  
However, using specific properties of the parameters $d=2$, $k=2$, $\ell=3$ Lovász and Yemini \cite{ly} \emph{do}
prove the Maxwell-Laman-theorem starting from an algebraic result in the same vein as our \lemref{ranknodrop},
providing evidence that our results may have some relevance to rigidity.

\end{document}